\documentclass[11pt,lefteqn]{article}

\usepackage{amsmath,amsfonts,amsthm,amssymb,graphics,epsfig,color,soul}
\usepackage{hyperref} 

\begin{document}
\date {\today}
\title{Optimal Polynomial Approximation to Rational Matrix Functions 
Using the Arnoldi Algorithm}

\author {Tyler Chen\footnote{NYU Courant, Department of Mathematics, 251 Mercer St., 
New York NY, 10003; NYU Tandon, Department of Computer Science and Engineering, 
370 Jay St. New York NY, 11201. 
email:  tyler.chen@nyu.edu}, Anne Greenbaum\footnote{University of Washington,
Applied Math Dept., Box 353925, Seattle, WA 98195.  email:  greenbau@uw.edu},
Natalie Wellen\footnote{University of Washington, 
Applied Math Dept., Box 353925, Seattle, WA 98195.  email:  nwellen@uw.edu}
}

\maketitle

\begin{abstract} 
Given an $n$ by $n$ matrix $A$ and an $n$-vector $b$, along with a rational function
$R(z) := D(z )^{-1} N(z)$, we show how to find the optimal approximation to $R(A) b$
from the Krylov space, $\mbox{span}( b, Ab, \ldots , A^{k-1} b)$, 
using the basis vectors
produced by the Arnoldi algorithm.  To find this optimal approximation requires
running $\max \{ \mbox{deg} (D) , \mbox{deg} (N) \} - 1$ extra Arnoldi steps and
solving a $k + \max \{ \mbox{deg} (D) , \mbox{deg} (N) \}$ by $k$ least squares problem.
Here {\em optimal} is taken to mean optimal in the 
$D(A )^{*} D(A)$-norm.  Similar to the case for linear systems, we show that 
eigenvalues alone cannot provide information about the convergence behavior of this 
algorithm and we discuss other possible error bounds for highly nonnormal matrices.
\end{abstract}

\paragraph{2020 Mathematical subject classifications\,:}65F10; 65F60; 47A12; 47A25

\noindent{\bf Keywords\,:}{ rational functions, Arnoldi algorithm}

\section{Introduction}
Given an $n$ by $n$ matrix $A$ and an $n$-vector $b$, along with a rational 
function $R(z) := D(z )^{-1} N(z)$, we wish to approximate $R(A)b$ using a polynomial
of degree at most $k-1$:  $R(A) b \approx P_{k-1} (A) b$.  A standard approach
is to use the Arnoldi algorithm to generate an orthonormal basis for the Krylov space 
$\mathcal{K}_k (A,b) := \mbox{span}( b, Ab, \ldots , A^{k-1} b )$. 
If $Q_k$ is the $n$ by $k$ matrix whose columns are these orthonormal basis vectors,
then the Arnoldi recurrence can be written in the form
\begin{equation}
A Q_k = Q_k H_k + h_{k+1,k} q_{k+1} e_k^T = Q_{k+1} H_{k+1,k} , \label{Arnoldi}
\end{equation}
where $e_k$ denotes the $k$th unit vector, $H_k$ is a $k$ by $k$
upper Hessenberg matrix, and $H_{k+1,k}$ is the $k+1$ by $k$ matrix obtained
by appending a row to $H_k$ consisting of $k-1$ zeros followed by the value
$h_{k+1,k}$.  

Given a function $f$, the standard way to approximate $f(A)b$
using vectors from the Krylov space is to take
\[
Q_k f( H_k ) Q_k^{*} b .
\]
This is called Arnoldi-FA (or Lanczos-FA when $A$ is Hermitian) in \cite{CGMM1}.
A standard result is that,
for any norm $\| \cdot \|$ induced by a matrix that commutes with $A$ and any polynomial $p$
of degree at most $k-1$, this approximation satisfies
\begin{eqnarray*}
\| f(A) b - Q_k f( H_k ) Q_k^{*} b \| & \leq &
\| f(A) b - p(A) b \| + \| p(A) b - Q_k p( H_k ) Q_k^{*} b \| + \| Q_k p( H_k ) Q_k^{*} b -
Q_k f( H_k ) Q_k^{*} b \| \\
 & = & \| f(A) b - p(A) b \| + 0 + \| Q_k ( p( H_k ) - f( H_k ) ) Q_k^{*} b \| \\
 & \leq & \left( \| f(A) - p(A) \|_2 + \| p( H_k ) - f( H_k ) \|_2 \right) \| b \| .
\end{eqnarray*}
Thus, the accuracy of the Arnoldi-FA approximation depends on how well the
{\em two} matrices $f(A)$ and $f( H_k )$ can be approximated by the same polynomial function
of the corresponding matrix.  Since the matrix $H_k$ is not known {\em a priori} and
depends on the right-hand side vector $b$, it would be preferable to have 
an estimate based only on how well $f(A)$ can be approximated by $p(A)$.
In particular, if the function $f$ is well-behaved around the eigenvalues of $A$,
but perhaps not at points in between, where eigenvalues of $H_k$ may lie,
the above error bound may be large and the error in the Arnoldi-FA 
approximation may be large.

In \cite{CGMM2}, an algorithm called Lanczos-OR was derived to generate the {\em optimal} 
(in a certain norm) approximation to $f(A)b$ from the Krylov space, when $A$ is a 
Hermitian matrix and $f$ is a rational function.  Here we extend this to non-Hermitian 
matrices, using a possibly simpler norm than the one used in \cite{CGMM2}.  
The algorithm, which we call Arnoldi-OR, is derived in Section 2.  

For diagonalizable matrices with fairly well-conditioned eigenvector matrices,
the error in the Arnoldi-OR approximation can be related to how well the rational function 
can be approximated by a polynomial of a certain degree $k-1$ at the eigenvalues of $A$.  
For nondiagonalizable or otherwise highly nonnormal matrices, however, this is not the case.  
In Section 3, we derive error bounds for highly nonnormal matrices based on the
numerical range $W(A) := \{ \langle Aq,q \rangle : \langle q,q \rangle = 1 \}$
and related sets in the complex plane, using results from \cite{CP} and \cite{CG,GW}.
Section 4 contains numerical examples.  Finally, in section 5, we show that,
as in the case of linear systems \cite{GPS} (i.e., for $R(A) = A^{-1}$), any 
nonincreasing convergence curve can be attained by the 2-norms of the Arnoldi-OR residuals 
for a matrix having any given eigenvalues.

\section{The Arnoldi-OR Algorithm}
The Arnoldi algorithm for generating an orthonormal basis for the Krylov space
$\mbox{span} ( b, Ab, A^2 b , \ldots )$ can be written as follows:

\vspace{.1in}
\begin{center}
\begin{tabular}{|l|} \hline
{\bf Arnoldi Algorithm:} \\
 \\
Set $q_1 = b / \| b \|_2$.  For $j=1, \ldots , k$, \\
~~~Set $\tilde{q}_{j+1} = A q_j$. \\
~~~For $i=1, \ldots , j$, \\
~~~~~~$h_{ij} = \langle \tilde{q}_{j+1} , q_i \rangle$; 
      $\tilde{q}_{j+1} \leftarrow \tilde{q}_{j+1} - h_{ij} q_i$; \\
~~~endFor \\
~~~Repeat orthogonalization if desired: \\
~~~For $i=1, \ldots , j$,  \\
~~~~~~$\tilde{q}_{j+1} \leftarrow \tilde{q}_{j+1} - \langle \tilde{q}_{j+1} ,
       q_i \rangle q_i$; \\
~~~endFor \\
~~~$h_{j+1,j} = \| \tilde{q}_{j+1} \|_2$; \\
~~~$q_{j+1} = \tilde{q}_{j+1} / h_{j+1,j}$. \\
endFor \\ \\
\hline
\end{tabular}
\end{center}
\vspace{.1in}

Forming the matrix $Q_{k+1}$ whose columns are the vectors $q_1 , \ldots , q_{k+1}$
and the upper Hessenberg matrix $H_{k+1,k}$ whose entries are the coefficients $h_{ij}$,
one obtains equation (\ref{Arnoldi}).  Assume wlog that $\| b \|_2 = 1$ so that 
$q_1 = b$.  Given the rational function
$R(z) := D(z )^{-1} N(z)$, we wish to find the vector $y$ that minimizes:
\begin{eqnarray}
\| R(A) b - Q_k y \|_{D(A )^{*} D(A)}^2 & = & \langle D(A )^{-1} N(A)b - Q_k y, D(A )^{*} D(A)
( D(A )^{-1} N(A) b - Q_k y ) \rangle \nonumber \\
 & = & \langle N(A)b - D(A) Q_k y , N(A) b - D(A) Q_k y \rangle \nonumber \\
 & = & \| N(A) Q_k e_1 - D(A) Q_k y \|_2^2 . \label{minimize1}
\end{eqnarray}
Note that if $D(A) = A$ and $N(A) = I$, this is the $2$-norm of the residual $b-A Q_k y$,
which is minimized by the GMRES algorithm.

Suppose $D(z) := \sum_{j=0}^J d_j z^j$ and $N(z) := \sum_{\ell =0}^L n_{\ell} z^{\ell}$.
Multiplying equation (\ref{Arnoldi}) by powers of $A$, we find
\begin{eqnarray}
A Q_k & = & Q_{k+1} H_{k+1,k} \nonumber \\
A^2 Q_k & = & A Q_{k+1} H_{k+1,k} = Q_{k+2} H_{k+2,k+1} H_{k+1,k} := Q_{k+2} H_{k+2,k} \nonumber \\
\vdots & & \nonumber \\
A^j Q_k & = & Q_{k+j} H_{k+j,k+j-1} H_{k+j-1,k+j-2} \ldots H_{k+1,k} := 
Q_{k+j} H_{k+j,k} . \label{AjQk}
\end{eqnarray}
Here we have defined the $k+j$ by $k$ matrices $H_{k+j,k} := H_{k+j,k+j-1} \ldots H_{k+1,k}$,
$j=1, 2, \ldots$, and we will define $H_{k,k} := I_{k \times k}$.
Substituting the expressions for $D(A)$ and $N(A)$ into (\ref{minimize1}), we can express the quantity
to be minimized as follows:
\begin{eqnarray*}
\| N(A) Q_k e_1 - D(A) Q_k y \|_2 & = &
\left\| \sum_{\ell = 0}^L n_{\ell} A^{\ell} Q_k e_1 - \sum_{j=0}^J d_j A^j Q_k y \right\|_2 \\
 & = & \left\| \sum_{\ell =0}^L n_{\ell} Q_{k+ \ell} H_{k+ \ell ,k} e_1 - ( \sum_{j=0}^J d_j Q_{k+j} H_{k+j,k} )
y \right\|_2 .
\end{eqnarray*}
Finally, let $\nu := \max \{ J,L \}$ and define $\hat{H}_{k+j,k}$ and $\hat{H}_{k+ \ell , k}$ to be the
$k+\nu$ by $k$ matrices obtained by appending rows of zeros to $H_{k+j,k}$ and $H_{k+ \ell ,k}$, resp.
Then we will choose $y$ to minimize
\[
\left\| Q_{k+ \nu} \left( \sum_{\ell =0}^L n_{\ell} \hat{H}_{k+ \ell , k} e_1 - ( \sum_{j=0}^J d_j 
\hat{H}_{k+j,k} ) y \right) \right\|_2 =
\left\| \left( \sum_{\ell =0}^L n_{\ell} \hat{H}_{k+ \ell , k} e_1 - ( \sum_{j=0}^J d_j 
\hat{H}_{k+j,k} ) y \right) \right\|_2 .
\]
This least squares problem can be solved using a standard QR factorization of the $k+ \nu$ by $k$
coefficient matrix, $\sum_{j=0}^J d_j \hat{H}_{k+j,k}$.

Note that $H_{k+j,k}$ is the top left $(k+j) \times k$ block of 
$H_{k+j}^j$, which is the same as the top left $(k+j) \times k$ block of $H_{k+ \nu}^j$
or of the $j$th power of any upper Hessenberg extension of $H_{k+j}$ such as the final 
upper Hessenberg matrix $H$ in the Arnoldi algorithm.
Also, the first $k$ entries of the remaining rows of $H_{k+ \nu}^j$ are zero, 
for $j \leq \nu$, so that $\hat{H}_{k+j,k}$ is the top left $(k+ \nu ) \times k$ block 
of $H_{k+ \nu}^j$.   
To see this, note that $H_{k+ \nu}$ can be written in the forms
\[
H_{k+ \nu} =
\left[ \begin{array}{cc} H_{k+ \nu , k+ \nu -1} & 
\begin{array}{c} \ast \\ \vdots \\ \ast \end{array} \end{array} \right] =
\left[ \begin{array}{cc} H_{k+ \nu -1, k+ \nu -2} &
\begin{array}{cc} \ast & \ast \\ \vdots & \vdots \\ \ast & \ast \end{array} \\
\begin{array}{ccc} 0 & \ldots & 0 \end{array} & 
\begin{array}{cc} \ast & \ast \end{array} \end{array} \right] = 
\left[ \begin{array}{cc} H_{k+ \nu -2, k+ \nu -3} &
\begin{array}{ccc} \ast & \ast & \ast \\ \vdots & \vdots & \vdots \\
\ast & \ast & \ast \end{array} \\
\begin{array}{ccc} 0 & \ldots & 0 \\ 0 & \ldots & 0 \end{array} &
\begin{array}{ccc} \ast & \ast & \ast \\ 0 & \ast & \ast \end{array} \end{array} \right]
= \ldots .
\]
Computing $H_{k+ \nu}^2$ by multiplying the first two representations, it is clear that
entries in the last column of the first representation are multiplied by zeros
so that the first $k+ \nu -2$ columns of $H_{k + \nu}^2$ are equal to
$H_{k+ \nu , k+ \nu -1} H_{k+ \nu -1 , k+ \nu -2} = H_{k+ \nu , k+ \nu -2}$.
Similarly, computing $H_{k+ \nu}^3$ by multiplying the three representations,
it is clear that the entries in the last two columns of the second representation
are multiplied by zeros so that the first $k+ \nu -3$ columns of $H_{k+ \nu}^3$
are equal to $H_{k+ \nu , k+ \nu -1} H_{k+ \nu -1, k+ \nu -2} H_{k+ \nu -2 , k+ \nu -3} =
H_{k+ \nu , k+ \nu - 3}$, etc. 

Therefore the least squares problem to be solved can be written in the form:
\[
\min_y \| N( H_{k+ \nu} )(:,1:k) e_1 - D( H_{k+ \nu} )(:,1:k) y \|_2 =
\min_y \| N( H_{k+ \nu} )(:,1) - D( H_{k+ \nu} )(:,1:k) y \|_2 ,
\]
where we have used MATLAB notation to indicate all rows and columns $1$ through $k$
of the matrix $D( H_{k+ \nu} )$ and the first column of $N( H_{k+ \nu} )$.

The basic Arnoldi-OR algorithm can be written as follows:

\vspace{.1in}
\begin{center}
\begin{tabular}{|l|} \hline
{\bf (Basic) Arnoldi-OR Algorithm} ($\nu = \max \{ \mbox{deg}(N), \mbox{deg}(D) \})$: \\
 \\
Run $1 + \nu$ steps of the Arnoldi algorithm to produce $Q_{1+ \nu}$ and $H_{1+ \nu}$. \\
Form the first column of $D( H_{1+ \nu} )$ and the first column of $N( H_{1+ \nu} )$. \\
~~Call these $\mathcal{D}$ and $\eta$, respectively. \\  
Solve the $(1+ \nu ) \times 1$ least squares problem $\mathcal{D} y \approx \eta$
for $y$, and use $\| \eta - \mathcal{D} y \|_2$ as the residual norm. \\
If the residual norm is sufficiently small, form the first approximation to $R(A) b$:
$x_1 = q_1 y$. \\
 \\
For $k=2: \mbox{kmax}$, \\
~~~Run step $k+ \nu$ of the Arnoldi algorithm to produce $Q_{k+ \nu}$ and $H_{k+ \nu}$.  \\
~~~Form the $k$th column and the $(k + \nu)$th row of $D( H_{k+ \nu} )(:,1:k)$. \\
~~~~~Append these as the last column and last row of $\mathcal{D}$.  
Append a $0$ to $\eta$. \\
~~~Solve the $(k+ \nu ) \times k$ least squares problem $\mathcal{D} y \approx \eta$ for $y$, 
and use $\| \eta - \mathcal{D} y \|_2$ as the residual norm. \\
~~~If the residual norm is sufficiently small, form the $k$th approximation to $R(A)b$: 
$x_k = Q_k y$. \\
endFor \\ \\
\hline
\end{tabular}
\end{center}
\vspace{.1in}

Although the matrix $\mathcal{D} = D( H_{k+ \nu} )(:,1:k)$ in the algorithm is
not upper Hessenberg if $D$ has degree greater than one, at each step, the new 
matrix $\mathcal{D}$ is obtained from the previous one by appending just one
row and one column.  If the least squares problems at previous steps have been
solved by applying Givens rotations to the coefficient matrix to make it 
upper triangular, applying the same Givens rotations to $\eta$, and then solving 
the resulting upper triangular linear system for $y$, then these same Givens rotations 
can be applied to the last column of $D( H_{k+ \nu} )(:,1:k)$, which is then
appended to the upper triangular matrix $\mathcal{D}$ from the previous step, 
along with a new row.  A new set of rotations must then be computed to eliminate 
entries in column $k$, rows $k+1 , \ldots , k+ \nu$.

Incorporating this reuse of rotations into the basic Arnoldi-OR algorithm,
we obtain the following algorithm:

\vspace{.1in}
\begin{center}
\begin{tabular}{|l|} \hline
{\bf Arnoldi-OR Algorithm} ($\nu = \max \{ \mbox{deg}(N), \mbox{deg}(D) \}$, 
$J = \mbox{deg}(D)$): \\
 \\
Run $1 + \nu$ steps of the Arnoldi algorithm to produce $Q_{1+ \nu}$ and $H_{1+ \nu}$. \\
Form the first column of $D( H_{1+ \nu} )$ and the first column of $N( H_{1+ \nu} )$. \\
Form Givens rotations $F_{1,2} , \ldots , F_{1,1+ J}$ to annihilate
entries $2$ through $1+ J$ in $D( H_{1+ \nu} )(:,1)$. \\
Set $\mathcal{D} := F_{1,2} \cdots F_{1,1+ J} D( H_{1+ \nu} )(:,1)$ \\
Apply these rotations to $N( H_{1+ \nu} )(:,1)$ to obtain
$\eta := F_{1,2} \cdots F_{1,1+ J} N( H_{1+ \nu} )(:,1)$. \\
Solve the upper triangular system $\mathcal{D}(1,1) y = \eta (1)$ for $y$,
and use $\| \eta - \mathcal{D} y \|_2$ as the residual norm. \\
If the residual norm is sufficiently small, form the first approximation to $R(A) b$:
$x_1 = q_1 y$. \\
 \\
For $k=2: \mbox{kmax}$, \\
~~~Run step $k+ \nu$ of the Arnoldi algorithm to produce $Q_{k+ \nu}$ and $H_{k+ \nu}$.  \\
~~~Form the $k$th column and the $(k + \nu )$th row of $D( H_{k+ \nu} )(:,1:k)$. \\
~~~Apply the previous rotations to the $k$th column of $D( H_{k+ \nu} )$: \\
~~~~~$( F_{k-1,k} \cdots F_{k-1,k-1+ J} ) \cdots ( F_{1,2} \cdots F_{1,1+ J} )
D( H_{k+ \nu} )(:,k)$. \\
~~~~~Append this as the last column of $\mathcal{D}$ and append the $(k+ \nu )$th row
of $D( H_{k+ \nu} )(:,1:k)$ \\
~~~~~as the last row of $\mathcal{D}$. Append a $0$ to $\eta$. \\
~~~Form Givens rotations $F_{k,k+1} , \ldots ,F_{k,k+ J}$ to annihilate entries
$k+1$ through $k+ J$ in column $k$ of $\mathcal{D}$. \\
~~~~~$\mathcal{D}(:,k) \leftarrow F_{k,k+1} \cdots F_{k,k+ J} \mathcal{D}(:,k)$. \\
~~~Apply these new rotations to $\eta$:  $\eta \leftarrow F_{k,k+1} \cdots F_{k,k+ J} 
\eta$. \\
~~~Solve the upper triangular system $\mathcal{D}(1:k,1:k) y = \eta (1:k)$ for $y$,
and use \\
~~~~~$\| \eta - \mathcal{D} y \|_2$ as the residual norm. \\
~~~If the residual norm is sufficiently small, form the $k$th approximation to $R(A)b$: 
$x_k = Q_k y$. \\
endFor \\ \\
\hline
\end{tabular}
\end{center}
\vspace{.1in}

Note that this assumes that $\mbox{kmax} + \nu \leq n$ (or the number of linearly 
independent vectors in the Krylov space ${\cal K}_n (A,b)$).  Once a full $n$ steps 
of the Arnoldi algorithm have been run and the final upper Hessenberg matrix $H$ produced,
the exact solution can be obtained by solving $D(H) y = N(H) e_1$ and setting
$x = Qy$.  A MATLAB code implementing the Arnoldi-OR algorithm can be found at:
\verb+ https://github.com/greenbau/Arnoldi-OR/tree/main +

\section{Error Bounds for Arnoldi-OR}
Let $S := D(A )^{*} D(A)$.
Since the Arnoldi-OR approximation $Q_k y$ is equal to $P_{k-1} (A) b$ where $P_{k-1}$ is the
$( k-1 )^{st}$ degree polynomial that minimizes the $S$-norm of the error,
we can write
\begin{equation}
\frac{\| R(A) b - Q_k y \|_2}{\| b \|_2} \leq \kappa ( S )^{1/2} \min_{p_{k-1} \in {\cal P}_{k-1}}
\| R(A) - p_{k-1} (A) \|_2 , \label{2normbound}
\end{equation}
where $\kappa ( S ) := \| S \|_2 \cdot \| S^{-1} \|_2$ is the condition number of $S$.
If $A$ is diagonalizable, with eigendecomposition $A = V \Lambda V^{-1}$, where
$\Lambda = \mbox{diag}( \lambda_1 , \ldots , \lambda_n )$, then
\begin{equation}
\frac{\| R(A) b - Q_k y \|_2}{\| b \|_2} \leq \kappa ( S )^{1/2} \kappa (V) \min_{p_{k-1} \in {\cal P}_{k-1}}
\max_{i=1, \ldots , n} | R( \lambda_i ) - p_{k-1} ( \lambda_i ) | . \label{2normevalbound}
\end{equation}

If $A$ is not diagonalizable or if $\kappa (V)$ is huge, 
one may still relate the error in the 
Arnoldi-OR approximation to how well one can approximate $R(z)$ on the {\em numerical range} of $A$,
$W(A) := \{ \langle Aq,q \rangle : \langle q,q \rangle = 1 \}$.  It is shown in \cite{CP} that
the numerical range is a $(1 + \sqrt{2} )$-spectral set, meaning that for any function $f$ that
is analytic in the interior of $W(A)$ and continuous on the boundary, 
$\| f(A) \| \leq (1 + \sqrt{2} ) \max_{z \in W(A)} | f(z) |$.  It follows
that the error in the Arnoldi-OR approximation satisfies
\begin{equation}
\frac{\| R(A) b - Q_k y \|_2}{\| b \|_2} \leq \kappa ( S )^{1/2} (1 + \sqrt{2} ) \min_{p_{k-1} \in {\cal P}_{k-1}}
\max_{z \in W(A)} | R( z ) - p_{k-1} ( z ) | . \label{2normWofAbound}
\end{equation}
If one prefers to work directly with the $S$-norm, the result in \cite{CP} still holds, provided the
$S$-inner product is used in the definition of the numerical range:  $W_{S} (A) := 
\{ \langle Aq, Sq \rangle : \langle q, Sq \rangle = 1 \} = 
\{ \langle S^{1/2} A S^{-1/2} ( S^{1/2} q ) , S^{1/2} q \rangle : \langle S^{1/2} q , S^{1/2} q \rangle = 1 \} =
W( S^{1/2} A S^{-1/2} )$.  This leads to the Arnoldi-OR error bound:
\begin{equation}
\frac{\| R(A) b - Q_k y \|_S}{\| b \|_S} \leq (1 + \sqrt{2} ) \min_{p_{k-1} \in {\cal P}_{k-1}}
\max_{z \in W( S^{1/2} A S^{-1/2} )} | R( z ) - p_{k-1} ( z ) | . \label{SnormWofAbound}
\end{equation}

The previous bounds are not useful if $R$ has a pole in $W(A)$ (or in $W( S^{1/2} A S^{-1/2} )$).
It is shown in \cite{CG} that one can remove a disk about such a pole $\xi \in W(A)$ of radius
$1/w(( \xi I - A )^{-1} )$, where $w( \cdot )$ denotes the {\em numerical radius} (the maximum
absolute value of all points in the numerical range), and still have a $(3 + 2 \sqrt{3})$-spectral set.
Replacing $W(A)$ in (\ref{2normWofAbound}) by $W(A)$ minus such a disk leads to the error bound: 
\begin{equation}
\frac{\| R(A) b - Q_k y \|_2}{\| b \|_2} \leq \kappa ( S )^{1/2} (3 + 2 \sqrt{3} ) \min_{p_{k-1} \in {\cal P}_{k-1}}
\max_{z \in W(A) \backslash \mathbb{D} ( \xi , 1/ w( ( \xi I - A )^{-1} ) )} 
| R( z ) - p_{k-1} ( z ) | . \label{2normWofAminusdiskbound}
\end{equation}
In \cite{GW}, this result is extended to allow the removal of $m$ disks, $\cup_{j=1}^m \mathbb{D} ( \xi_j ,
1/w( ( \xi_j I - A )^{-1} ) )$, with the subsequent error bound
\[
\frac{\| R(A) b - Q_k y \|_2}{\| b \|_2} \leq \kappa ( S )^{1/2} \left( 1+2m + \sqrt{(1+2m )^2 + 2m+1}  \right) \cdot 
\]
\begin{equation}
\min_{p_{k-1} \in {\cal P}_{k-1}}
\max_{z \in W(A) \backslash \cup_{j=1}^m \mathbb{D} ( \xi_j , 1/ w( ( \xi_j I - A )^{-1} ) )} 
| R( z ) - p_{k-1} ( z ) | . \label{2normWofAminusmdisksbound}
\end{equation}

Note that each of the error bounds (\ref{2normbound} - \ref{2normWofAminusdiskbound}) 
involves information beyond just the eigenvalues of $A$.  We will show in section
\ref{sec:any-convergence} that eigenvalues alone cannot provide good error bounds;
as is the case with linear systems, any nonincreasing convergence curve can be
obtained with Arnoldi-OR for a matrix having any given eigenvalues.

\section{Numerical Examples} \label{sec:examples}
Here we present numerical experiments, comparing Arnoldi-OR with Arnoldi-FA and with
the optimal (in 2-norm) approximation from the Krylov space (which we compute
by directly evaluating $R(A)b$ and finding its orthogonal projection onto the 
Krylov space).  We compare the 2-norm of the residual $\| N(A) b - D(A) x_k \|_2$
and the 2-norm of the error $\| D(A )^{-1} N(A) b - x_k \|_2$ for each of these
approximations.  We also consider various error bounds for these quantities.
For a comprehensive comparison of different error bounds for GMRES, see
\cite{Embree}.

\subsection{Random $A$} \label{subsec:random}
Let $A$ be a $100$ by $100$ real random matrix (entries from a standard normal
distribution) plus a multiple of the identity, and let $b$ be a real random vector 
of length $100$, which we normalize for convenience.  Here we take $D$ to be a random 
cubic and $N$ to be a random quadratic:
\[
D(z) = \gamma (z - r_1 )(z - r_2 )(z - r_3 ),~~
N(z) = \delta (z - s_1 )(z - s_2 ) ,
\]
where $\gamma , r_1 , r_2 , r_3 , \delta , s_1 , s_2$ are complex numbers with
real and imaginary parts chosen from a standard normal distribution.
Depending on the multiple of the identity that is added to $A$, the roots of $D$
may lie inside or outside the numerical range of $A$.

\begin{figure}[ht]
\centerline{\epsfig{file=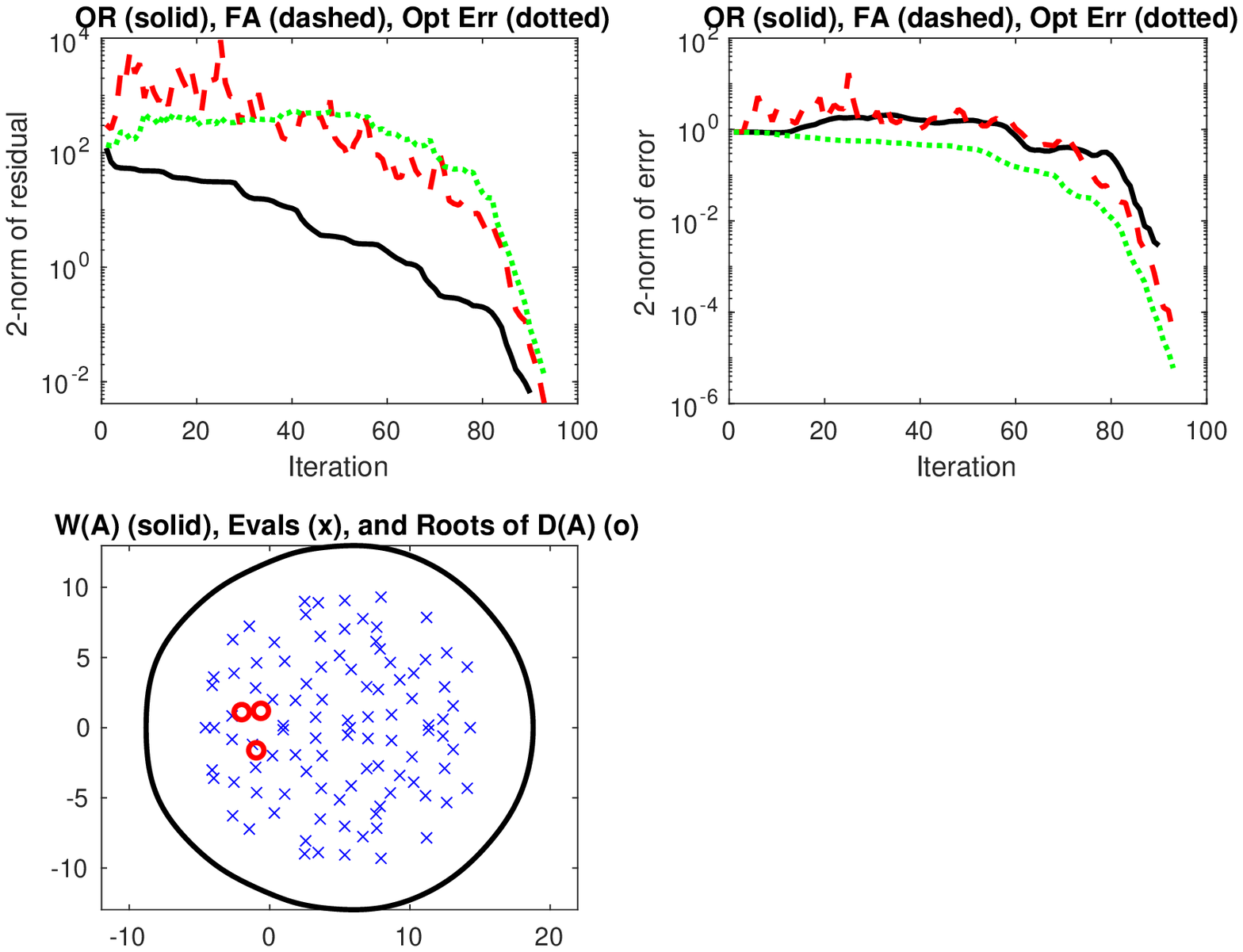,width=3.6in}}
\caption{$A = \mbox{randn}(100) + 5 \mbox{I}$, $b$ random, $D(z)$ random cubic,
$N(z)$ random quadratic.  Top plots show $2$-norm of residual and $2$-norm of error 
in Arnolid-OR (solid, black), Arnoldi-FA (dashed, red), and optimal approximation
in 2-norm (dotted, green).  Lower left plot shows eigenvalues (x, blue) and numerical
range (solid, black boundary curve) of $A$ and poles of $R(z)$ (o, red).
\label{fig:test1afig}}
\end{figure}

Figure \ref{fig:test1afig} shows the results of applying Arnoldi-OR (solid, black) and
Arnoldi-FA (dashed, red) to approximate $D(A )^{-1} N(A) b$ when $A = \mbox{randn}(100) + 
5 \mbox{I}$.  The condition number of a matrix $V$ of unit eigenvectors of $A$ was $98$,
and the condition number of $D(A)$ was about $2.3e+3$.  The top left plot
shows the 2-norm of the residual $\| N(A) b - D(A) x_k \|_2$ at each step, and the
top right plot shows the 2-norm of the error $\| D(A )^{-1} N(A) b - x_k \|_2$.
Also shown is the 2-norm of the residual and error in the orthogonal projection
of the true solution onto the Krylov subspace (dotted, green).  The 2-norm of the
residual in the Arnoldi-OR approximation is (necessarily) the smallest, while 
the 2-norm of the error in the orthogonal projection of the true solution is (necessarily)
the smallest.
One can see that Arnoldi-FA does not converge monotonically in either of these norms.
The lower left plot shows the eigenvalues (blue x's) and the numerical range (solid,
black boundary curve) of $A$, and it also shows the roots of $D(z)$ (i.e., the poles 
of $R(z) := D(z )^{-1} N(z)$) (red o's).  The poles of the rational function lie within the 
numerical range and close to some of the eigenvalues, and convergence is slow.
One could use inequality (\ref{2normWofAminusmdisksbound}) to obtain an upper bound
on the 2-norm of the error at each step, but computing the best uniform polynomial
approximation to $R(z)$
on the multiply connected domain consisting of $W(A)$ minus a disk about each pole
is difficult.  Since the eigenvector matrix of $A$ is not very ill-conditioned,
bound(\ref{2normevalbound}) is likely to be a more realistic estimate.

\begin{figure}[ht]
\centerline{\epsfig{file=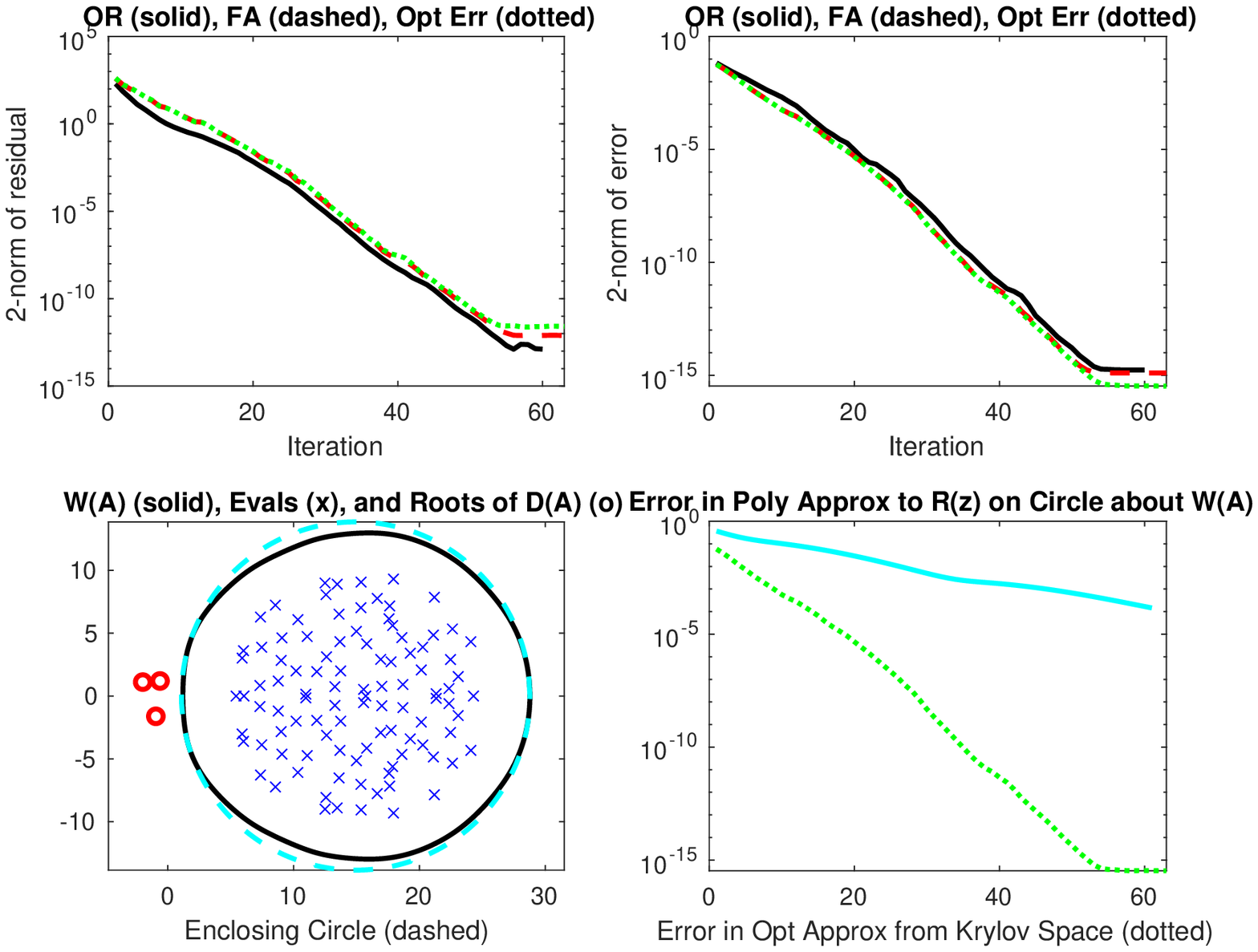,width=3.6in}}
\caption{$A = \mbox{randn}(100) + 15 \mbox{I}$, $b$ random, $D(z)$ random cubic,
$N(z)$ random quadratic.  Top plots show $2$-norm of residual and $2$-norm of error 
in Arnolid-OR (solid, black), Arnoldi-FA (dashed, red), and optimal approximation
in 2-norm (dotted, green).  Lower left plot shows eigenvalues (x, blue) and numerical
range (solid, black boundary curve) of $A$ and poles of $R(z)$ (o, red).  It also
shows the smallest circle enclosing $W(A)$ (dashed, cyan).
Lower right plot shows error in near best polynomial approximation to $R(z)$
on smallest circle enclosing $W(A)$ (solid, cyan), and how it compares to error
in optimal approximation to $R(A)b$ from successive Krylov spaces (dotted, green).
\label{fig:test1bfig}}
\end{figure}

Figure \ref{fig:test1bfig} shows the corresponding results when $A$ has the same 
random part but we now add $15 \mbox{I}$.  Here $b$, $D$, and $N$ are the same as 
in Figure \ref{fig:test1afig}.  Convergence is faster and there is less difference
between the convergence curves of Arnoldi-OR, Arnoldi-FA, and the optimal approximation
in the 2-norm, as shown in the top plots.
The condition number of the eigenvector matrix of $A$ is, of course, the same, but
now the condition number of $D(A)$ is $228$.
Now the poles of $R$ lie outside $W(A)$, as seen in the lower left plot.  Inequality
(\ref{2normWofAbound}) now gives an upper bound on the 2-norm of the error in terms of
the error in the best polynomial approximation to $R(z)$ on $W(A)$.  
The best uniform polynomial approximation on a set in the
complex plane is difficult to compute, but a near best approximation can be obtained
from the Faber series for $R(z)$.  See, for example, \cite[Ch.~1, Sec.~6]{Gaier}.  
The Faber series requires a conformal
map from the exterior of the unit disk to the exterior of $W(A)$.  This can be computed
using the Schwarz-Christoffel package of Driscoll \cite{Driscoll}.  For this problem,
however, since $W(A)$ is very close to a disk, we enclosed $W(A)$ in a disk (dashed,
cyan curve in the lower left plot), and used the Faber polynomials for the disk, which
are very simple.  Once the disk is mapped to the unit disk,
the Faber polynomials are just the powers
of $z$, and the Faber series for $R(z)$ is just the Taylor series.  The solid cyan curve
in the lower right plot shows
\[
\max_{z \in \mathbb{D}(c,r)} | R(z) - p_{k-1} (z) | ,
\]
where $\mathbb{D}(c,r)$, which denotes the disk with center $c$ and radius $r$,
is taken to be the smallest disk containing $W(A)$, and
$p_{k-1} (z)$ is the sum of the first $k$ terms in the Faber series for $R(z)$
on $\mathbb{D}(c,r)$.  Comparing this to the dotted green curve showing the error in the optimal
2-norm approximation to $R(A)b$ from successive Krylov spaces, one sees that there is a large
difference.  The reason is that the numerical range (or the circle enclosing the numerical
range) comes {\em much} closer to the poles of $R$ than do the eigenvalues of $A$.  If one
instead computes the best uniform polynomial approximation to $R$ on the smallest 
disk enclosing 
the spectrum of $A$, one sees a curve that is much closer to the dotted green curve. 
Since the eigenvectors of $A$ are not very ill-conditioned, inequality (\ref{2normevalbound})
probably provides the best error bound for this problem.

\begin{figure}[ht]
\centerline{\epsfig{file=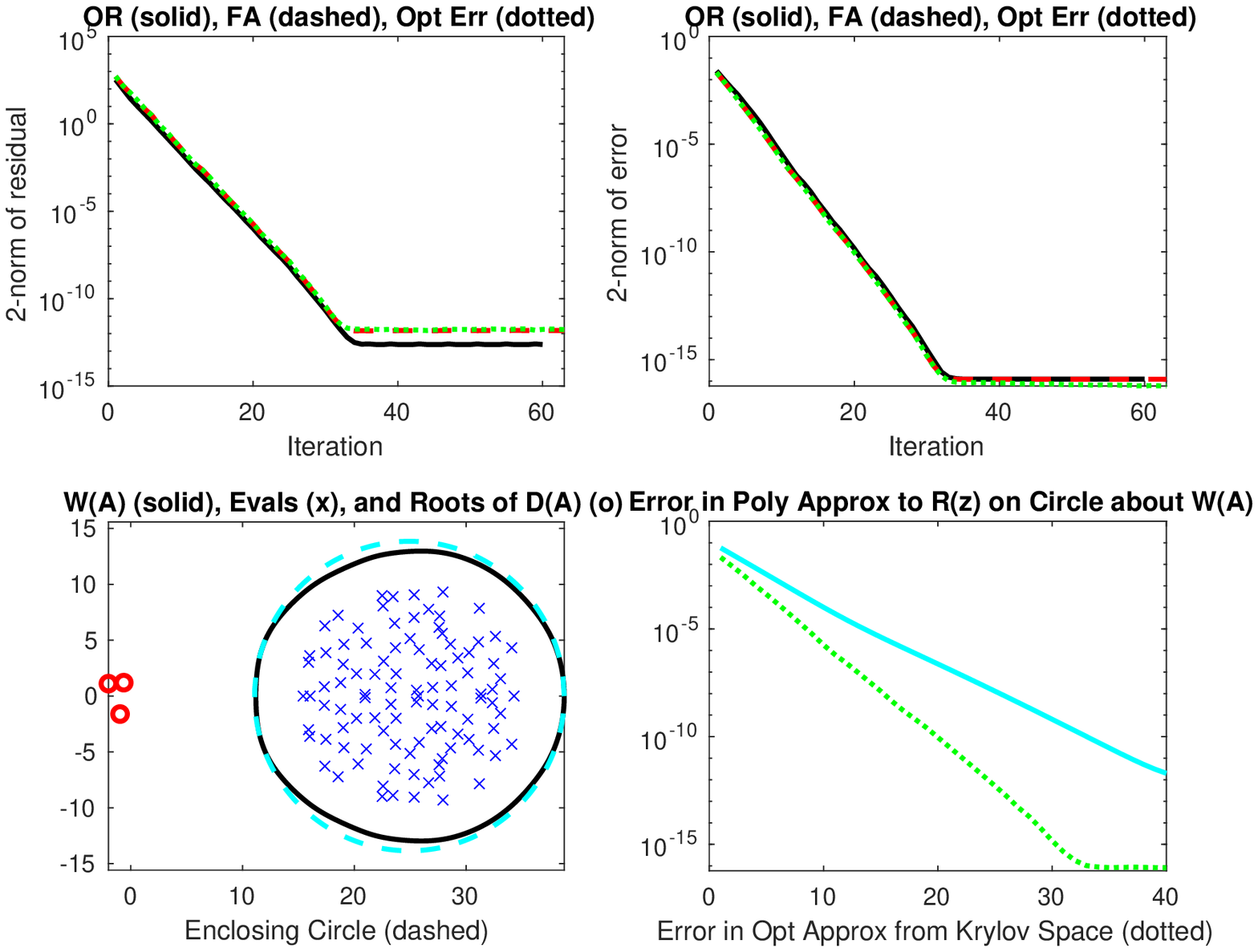,width=3.6in}}
\caption{$A = \mbox{randn}(100) + 25 \mbox{I}$, $b$ random, $D(z)$ random cubic,
$N(z)$ random quadratic.  Top plots show $2$-norm of residual and $2$-norm of error 
in Arnolid-OR (solid, black), Arnoldi-FA (dashed, red), and optimal approximation
in 2-norm (dotted, green).  Lower left plot shows eigenvalues (x, blue) and numerical
range (solid, black boundary curve) of $A$ and poles of $R(z)$ (o, red).  It also
shows the smallest circle enclosing $W(A)$ (dashed, cyan).
Lower right plot shows error in near best polynomial approximation to $R(z)$
on smallest circle enclosing $W(A)$ (solid, cyan), and how it compares to error
in optimal approximation to $R(A)b$ from successive Krylov spaces (dotted, green).
\label{fig:test1cfig}}
\end{figure}

Figure \ref{fig:test1cfig} shows results when $A$ has the same random part but
we now add $25 \mbox{I}$.  The condition number of $D(A)$ is $24.6$.
Now convergence is significantly faster and all of the
convergence curves in the upper plots are very close to each other.  Now the poles of $R$
lie further outside $W(A)$, as seen in the lower left plot.  A disk containing $W(A)$
(dashed, cyan curve in lower left plot) does not come so much closer, in a relative
sense, to these poles, than the eigenvalues themselves do.  The result in the lower
right plot is that the error bound based on the truncated Faber series approximation
to $R(z)$ on this disk is not so much larger than the actual $2$-norm of the error
in the best approximation to $R(A)b$ from successive Krylov spaces.

\subsection{Partial Fractions Decomposition}
A rational function can always be decomposed in a partial fractions decomposition.
For example, if the degree of $N$ is less than the degree of $D$ and the roots
of $D$ are all simple, then
\[
\frac{N(z)}{D(z)} = \sum_{i=1}^{\mbox{deg}(D)}
\frac{N( r_i ) / D' ( r_i )}{z - r_i} ,
\]
where the $r_i$'s are the roots of $D$.  Instead of applying Arnoldi-OR
directly to the function $R(z) = N(z)/D(z)$, one might apply it to
each term $(N( r_i ) / D' ( r_i )) /(z - r_i )$ and add the results.
Since each term $(N( r_i ) / D' ( r_i )) (A - r_i I )^{-1} b$
is approximated by a linear combination of the same Krylov vectors, 
$\{ b , Ab , \ldots , A^{k-1} b \}$,
an orthogonal basis for this space need be generated only once and the approximations
to different terms can be computed simultaneously, using GMRES (i.e., Arnoldi-OR for 
linear systems.) See, for example \cite{Napoli}.  For the
2-norm of the residual, $\| N(A) b - D(A) x_k \|_2$, Arnoldi-OR applied directly
to $N(z)/D(z)$ generates the optimal approximation from the Krylov space
and so cannot be beat in terms of number of iterations.  If $D(A)$ is
ill-conditioned, however, then there may be a large difference between the
2-norm of the residual and, say, the 2-norm of the error, 
$\| D(A )^{-1} N(A) b - x_k \|_2$, and using the partial fractions decomposition
might prove advantageous.

\begin{figure}[ht]
\centerline{\epsfig{file=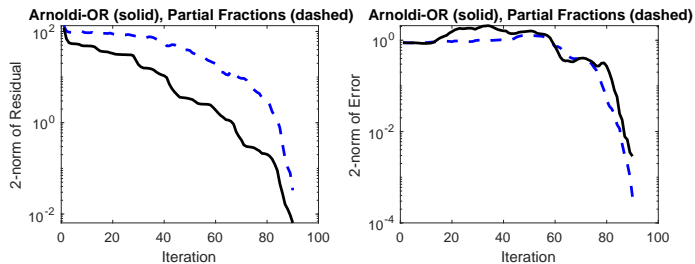,width=3.6in}}
\caption{$A = \mbox{randn}(100) + 5 \mbox{I}$, $b$ random, $D(z)$ random cubic, 
$N(z)$ random quadratic.
Left plot shows $2$-norm of the residual and right plot shows the $2$-norm of
the error in Arnoldi-OR applied to $N(z)/D(z)$ (solid, black) and in Arnoldi-OR 
applied to each term in the partial fractions decomposition (dashed, blue).  
\label{fig:pfrac}}
\end{figure}

Figure \ref{fig:pfrac} shows a comparison of Arnoldi-OR applied directly to
the first problem of the previous subsection ($A$ a random $100 \times 100$ matrix
plus $5$ times the identity, $N$ a random quadratic, and $D$ a random cubic) and
Arnoldi-OR applied to each term in the partial fractions decomposition of
$N(z)/D(z)$.  As noted above, the 2-norm of the residual in Arnoldi-OR 
applied directly to $N(z)/D(z)$ is always smaller than (or at least as small as) that 
in the computation based on partial fractions, but in this example, the 2-norm
of the error in the two approximations is very similar, with that in the
partial fractions computation sometimes being slightly smaller.

\subsection{Grcar Matrix}
The Grcar matrix is a nonsymmetric Toeplitz matrix with $-1$'s on the subdiagonal
and $1$'s on the main diagonal and the next three superdiagonals.  It can be
generated in MATLAB by typing \verb+gallery('grcar',n)+.  It has very ill-conditioned
eigenvalues, so that the error bound (\ref{2normevalbound}) is not useful.

\begin{figure}[ht]
\centerline{\epsfig{file=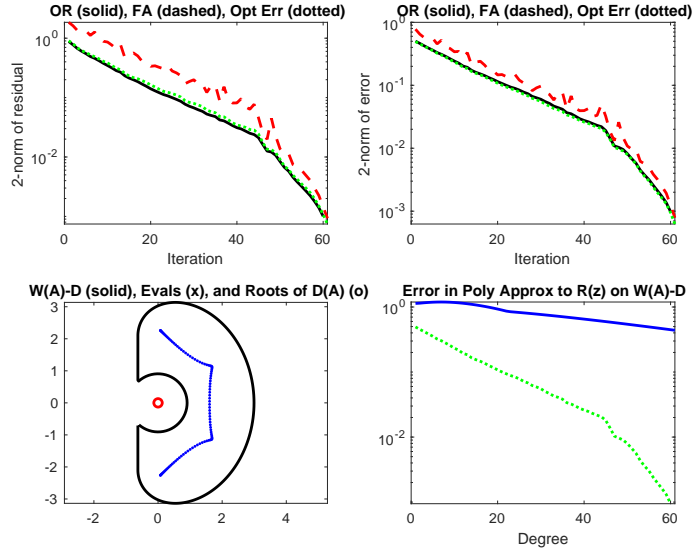,width=3.6in}}
\caption{$A = \mbox{gallery('grcar',100)}$, $b$ random, $D(z) = z$, $N(z) = 1$.
Top plots show $2$-norm of residual and $2$-norm of error 
in Arnolid-OR (solid, black), Arnoldi-FA (dashed, red), and optimal approximation
in 2-norm (dotted, green).  Lower left plot shows eigenvalues (dots, blue) and
$W(A) \backslash \mathbb{D}(0, 1/ w( A^{-1} ))$ (solid, black boundary curve).
Pole of $R(z)$ (o, red) is at the origin.  
Lower right plot shows error in near best polynomial approximation to $R(z)$
on $W(A) \backslash \mathbb{D}(0, 1/ w( A^{-1} ))$ (solid, blue), along with error
in optimal approximation to $R(A)b$ from successive Krylov spaces (dotted, green).
\label{fig:test2fig}}
\end{figure}

Here we considered only the problem of solving a linear system:  $D(A) = A$,
$N(A) = I$.  The upper plots in Figure \ref{fig:test2fig} show the convergence of 
Arnoldi-OR (solid, black), Arnoldi-FA (dashed, red), and the optimal approximation 
to $A^{-1} b$ in the 2-norm (dotted, green), where
$A$ is a $100$ by $100$ Grcar matrix and $b$ is a real random (normalized) vector
of length $100$.  The lower left plot shows the eigenvalues (blue dots) and
$W(A) \backslash \mathbb{D}(0, 1/ w( A^{-1} ))$ (solid, black boundary curve).
The pole of $R(z)$ (i.e., the origin) lies inside $W(A)$, so error bound 
(\ref{2normWofAbound}) is not useful, but error bound (\ref{2normWofAminusdiskbound})
can be applied.  We used the \verb+sc+ package \cite{Driscoll} to compute the
Faber polynomials for the region pictured in the lower left plot and from these
we derived the Faber series for $f(z) = 1/z$.  The error in approximating 
$f(z)$ by the first $k$ terms in this series is plotted in the lower right
(solid, blue), along with the error in the optimal 2-norm approximation to
$A^{-1} b$ from successive Krylov spaces (dotted, green).  The error bound
based on Faber polynomials is a large overestimate, but at least it does 
indicate convergence, while none of the other error bounds would be less than one.

\subsection{Matrix Exponential}
It was noted in subsection \ref{subsec:random} that when the poles of the rational
function lie well away from the eigenvalues and the numerical range of $A$, 
there may be little difference between the behavior of Arnoldi-OR, Arnoldi-FA,
and the optimal 2-norm approximation from the Krylov space.  Many functions, 
such as the matrix exponential $\exp (A)$ can be well-approximated by
rational functions whose poles may lie well outside the numerical range of $A$.

\begin{figure}[ht]
\centerline{\epsfig{file=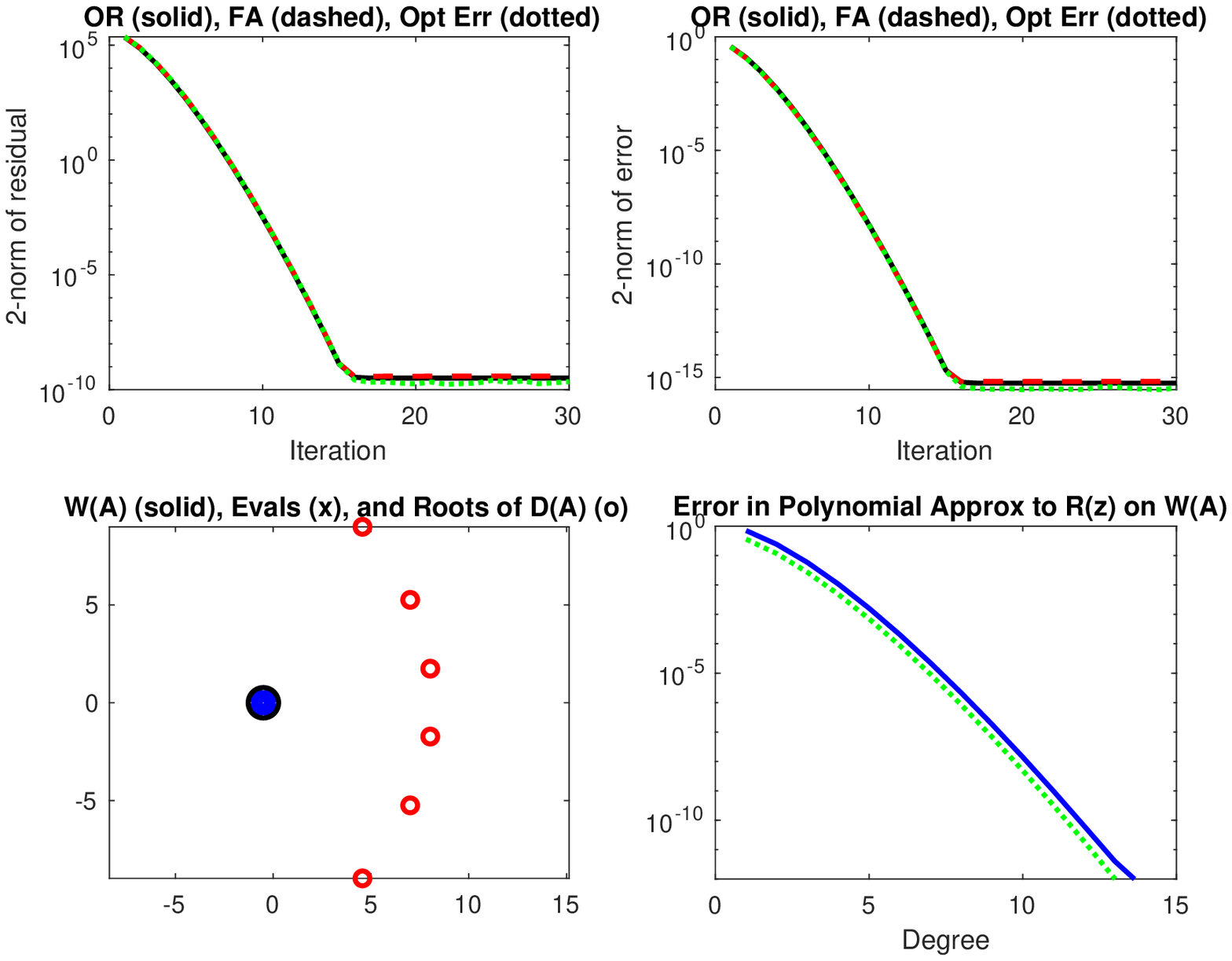,width=3.6in}}
\caption{$A = \mbox{transient\_demo(50)}$, $b$ random, $D$ and $N$ of degree 6.
Top plots show $2$-norm of residual and $2$-norm of error 
in Arnolid-OR (solid, black), Arnoldi-FA (dashed, red), and optimal approximation
in 2-norm (dotted, green) to rational function approximation of $\exp (A)b$.  
Lower left plot shows eigenvalues (x, blue) and
$W(A)$ (solid, black boundary curve), and poles of $R(z)$ (o, red).
Lower right plot shows error in best polynomial approximation to $R(z)$
on $W(A)$ (solid, blue) along with error in ortimal polynomial approximation
to $R(A)b$ (dotted, green).
\label{fig:transientexp}}
\end{figure}

Here we used routine \verb+rkfit+ from the \verb+RKToolkit+ \cite{RKToolkit} to find
a good rational function approximation to $\exp (A) b$, where $A$ was the
\verb+transient_demo+ matrix from \verb+eigtool+ \cite{eigtool}.  We set the degree 
of the numerator and denominator to $6$ and obtained a rational approximation that 
matched $\exp (A) b$ to near the machine precision.  This matrix is highly nonnormal
and demonstrates how $\| \exp (tA) \|_2$ can grow with $t > 0$ before decaying to $0$,
even though the eigenvalues lie in the left half-plane.  With the matrix size $n=50$,
we obtained the results in Figure \ref{fig:transientexp}.  The top plots show the
convergence of Arnoldi-OR (solid, black), Arnoldi-FA (dashed, red), and the optimal
approximation in 2-norm to the rational function approximation of $\exp (A) b$ 
(dotted, green).  In this case the convergence
curves are all very close in terms of both the 2-norm of the residual and the 2-norm
of the error.  The condition number of $D(A)$ was $2.17$.  The bottom left plot shows
the eigenvalues (x, blue), the numerical range (solid black boundary curve), and the
poles of the rational function (o, red).  The bottom right plot 
shows the error in the best polynomial approximation to $R(z)$ on $W(A)$ (solid, blue)
along with the error in the optimal polynomial approximation to $R(A)b$
(dotted, green).

While all methods converged very well when computing $\exp (A) b$, it can be 
much more challenging to compute $\exp (tA) b$ for $t >> 1$.  The 2-norm
of $\exp (tA )$ grows to about $4000$ around $t=60$, before decreasing towards $0$.
If we try to approximate $\exp (tA) b$ for larger 
values of $t$, then it is more difficult to find a good rational function 
approximation to $\exp (tA) b$.  The rational function must have higher degree
and more ill-conditioned numerator and denominator.
The problem of approximating this rational
function of $A$ with a polynomial in $A$ also becomes more difficult.
For some alternative ways to approximate $\exp (tA)$ or $\exp (tA)b$ for larger 
values of $t$, see, for example, \cite{Moler}.

\begin{figure}[ht]
\centerline{\epsfig{file=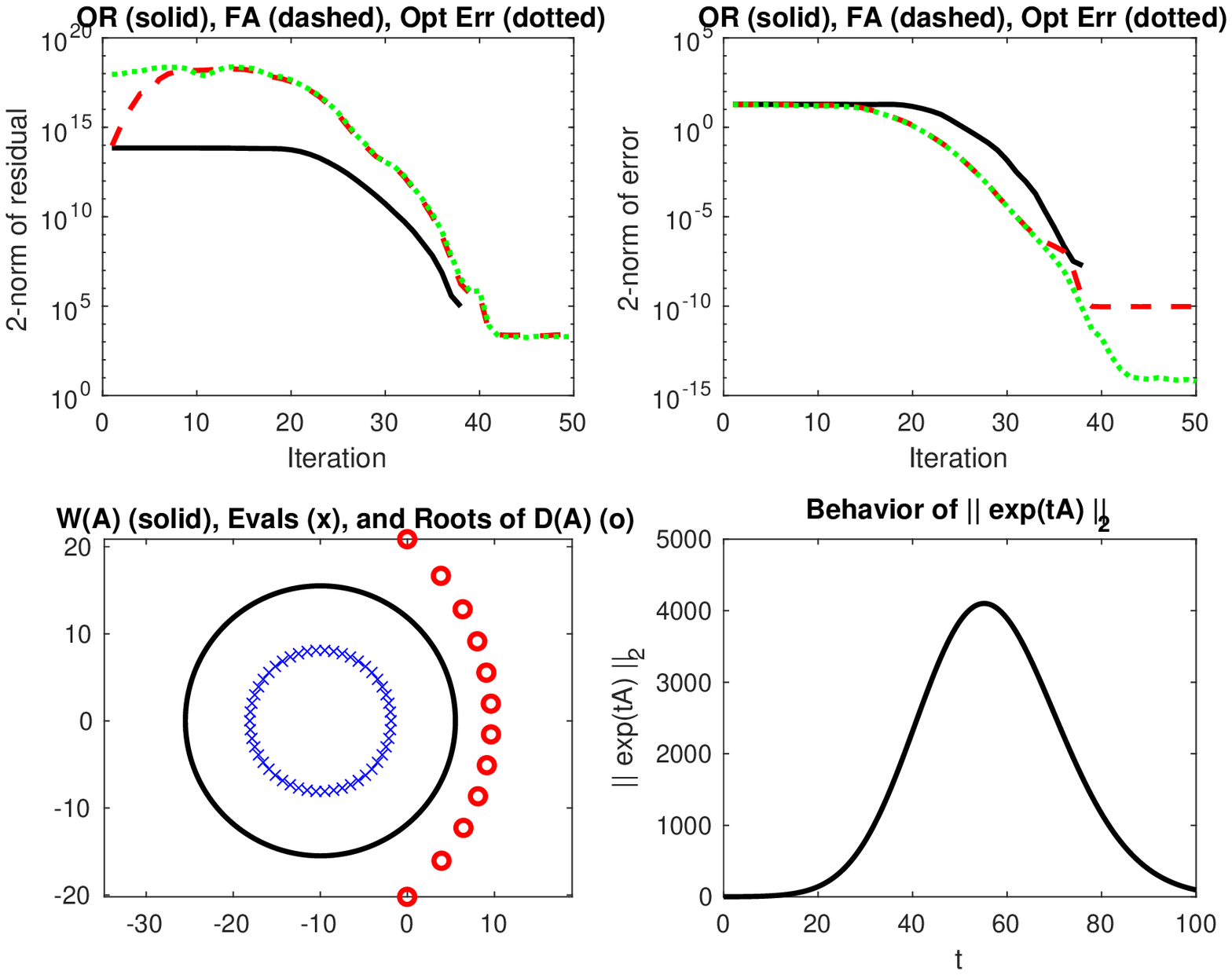,width=3.6in}}
\caption{$A = \mbox{transient\_demo(50)}$, $b$ random, $D$ and $N$ of degree 12.
Top plots show $2$-norm of residual and $2$-norm of error 
in Arnolid-OR (solid, black), Arnoldi-FA (dashed, red), and optimal approximation
in 2-norm (dotted, green) to rational function approximation of $\exp (20A)$.  
Lower left plot shows eigenvalues of $20A$ (x, blue) and
$W(20 A)$ (solid, black boundary curve), and poles of $R(z)$ (o, red).
Lower right plot shows behavior of $\| \exp (tA) \|_2$, $0 \leq t \leq 100$.
\label{fig:transientexp2}}
\end{figure}

Figure \ref{fig:transientexp2} shows the convergence of Arnoldi-OR, Arnoldi-FA,
and the optimal 2-norm approximation to a rational function approximation of
$\exp (20A) b$.  The rational function had numerator and denominator of degree $12$
and differed from $\exp (20A) b$ by about $1.9e-8$.  The denominator $D(A)$ had
condition number about $1.6e+6$, and the norm of $N(A) b$ was about $7.0e+13$.  
While Arnoldi-OR produced the smallest residual
norm, $\| N(A) b - D(A) x_k \|_2$, at each iteration (which, in an absolute
sense, was still quite large), it was the {\em least} effective
in reducing the $2$-norm of the error, $\| D(A )^{-1} N(A) b - x_k \|_2$.  
In this case
Arnoldi-FA produced results much closer to the optimal 2-norm approximation to
the rational function from the Krylov space, as illustrated in the top right
plot of Figure \ref{fig:transientexp2}.  In the bottom left plot, we see that
the poles of the rational function approximation are now closer to the eigenvalues
and numerical range of the matrix $20A$.  The bottom right plot shows the growth
of $\| \exp (tA) \|_2$ for $t$ running from $0$ to $100$.

\section{Constructing a problem with a given convergence curve and any prescribed 
eigenvalues} \label{sec:any-convergence}

In this section, we will assume that $N(A) = I$.  

Let $\{ \lambda_1 , \ldots , \lambda_n \}$ be a set of points in the complex
plane, with $D( \lambda_j ) \neq 0$, $j=1, \ldots , n$, and let $\varphi (0) \geq 
\varphi (1) \geq \ldots \geq \varphi (n- J ) > 0$ be any nonincreasing sequence.  
We will construct a matrix $A$ with eigenvalues $\lambda_1 , \ldots , \lambda_n$ 
and a right hand side vector $b$ for which the residual norm, 
$\| b - D(A) x_k \|_2$, is equal to $\varphi (k)$, $k=0,1 , \ldots , n - J$.
Our approach follows closely \cite[Section 2]{GPS}.

Note that the Arnoldi-OR residual $r_k := b - D(A) x_k$ satisfies the minimization 
condition
\begin{equation}
\| r_k \|_2 = \min_{u \in D(A) \mathcal{K}_k (A,b)} \| b - u \|_2 .
    \label{eqn:res_cond}
\end{equation}
Let $v_1 , \ldots , v_k$ be an orthonormal basis for $D(A) \mathcal{K}_k (A,b)$,
$k=1, \ldots , n-J$, and let $v_{n- J +1} , \ldots , v_n$ be any orthonormal
completion of $v_1 , \ldots , v_{n- J}$.  
Then (\ref{eqn:res_cond}) implies that 
\[
u = \sum_{j=1}^k \langle b , v_j \rangle v_j ,~~~
r_k = \sum_{j=k+1}^{n} \langle b , v_j \rangle v_j ,
\]
and
\[
\| r_k \|_2^2 = \sum_{j=k+1}^{n} | \langle b , v_j \rangle |^2 .
\]
It follows that
\[
| \langle b , v_k \rangle | = \sqrt{\| r_{k-1} \|_2^2 - \| r_{k} \|_2^2} ,~~ k=1, \ldots ,
n- J .
\]

From the given sequence $\varphi (0) , \ldots , \varphi (n- J )$, 
define the differences
\[
\psi (k) := \sqrt{ \varphi (k-1 )^2 - \varphi (k )^2} ,~~k=1, \ldots , n- J .
\]
We will construct $A$ with eigenvalues $\lambda_1 , \ldots , \lambda_n$,
$b$ with $\| b \|_2 = \varphi (0)$ and orthonormal vectors $v_1 , \ldots , v_{n-J}$,  
so that 
\begin{equation}
\mbox{span} \{ v_1 , \ldots , v_k \} = D(A) \mathcal{K}_k (A,b) ,~~k=1, \ldots , n-J
\label{spancond}
\end{equation}
and
\begin{equation}
| \langle b , v_k \rangle | = \psi (k),~~ k=1, \ldots , n-J . \label{residnorms}
\end{equation}

Let $v_1 , \ldots , v_{n-J}$ be any orthonormal vectors in $\mathcal{C}^n$ and
choose $A$ to satisfy 
\begin{equation}
A v_i = v_{i+1},~~ i=1, \ldots , n-J-1 \label{firstcond}
\end{equation}
Choose $b$ to satisfy (\ref{residnorms}) and the condition $\| b \|_2 = \varphi (0)$.
This means that $b$ is linearly independent from the vectors $v_1 , \ldots , v_{n-J}$,
since $\| [ v_1 , \ldots , v_{n-J} ]^{*} b \|_2^2 = \varphi(0 )^2 - \varphi (n-J )^2$
and $\varphi (n-J) > 0$.  If $J=1$, then $b$ is uniquely determined (up to multiplication
by a unit scalar) by these conditions, but if $J > 1$, then any vector $b$ satisfying
these conditions can be used.

If $J = 1$, so that $D(z) = c(z - \gamma )$, choose $A$ to satisfy, in addition to 
(\ref{firstcond}),
\begin{equation}
A b = v_1 + \gamma b , \label{secondcond}
\end{equation}
so that $v_1 = c^{-1} D(A) b$.  Then (\ref{firstcond}) implies that (\ref{spancond})
holds.
 
If $J > 1$, having chosen the orthonormal vectors $v_1 , \ldots , v_{n-J}$
and a vector $b$ satisfying (\ref{residnorms}), 
let $u_1 , \ldots , u_{J-1}$ be any vectors that are linearly
independent of $\{ v_1 , \ldots , v_{n-J} , b \}$.
Factor $D(z) = c (z - \gamma_1 ) \cdots (z - \gamma_J )$, and
choose $A$ to satisfy, in addition to (\ref{firstcond}),
\begin{eqnarray}
A b & = & u_1 + \gamma_1 b \nonumber \\
A u_1 & = & u_2 + \gamma_2 u_1 \nonumber \\
 & \vdots & \nonumber \\
A u_{J-2} & = & u_{J-1} + \gamma_{J-1} u_{J-2} \nonumber \\
A u_{J-1} & = & v_1 + \gamma_J u_{J-1} . \label{secondcondb}
\end{eqnarray}
Then $v_1 = (A - \gamma_J I) u_{J-1} = \ldots = \prod_{j=1}^J (A - \gamma_j I ) b = 
c^{-1} D(A)b$, and, again, (\ref{firstcond}) implies that (\ref{spancond}) holds.

We have now defined the action of $A$ on $n-1$ vectors via (\ref{firstcond})
and (\ref{secondcond}) or (\ref{secondcondb}).  We complete the definition of $A$
by defining
\begin{equation}
A v_{n-J} = \beta_0 b + \sum_{j=1}^{J-1} \beta_j u_j + 
\sum_{j=J}^{n-1} \beta_j v_{j-J+1} . \label{thirdcond}
\end{equation}

In the basis $\mathcal{B} := \{ b, u_1 , \ldots , u_{J-1} , v_1 , \ldots , v_{n-J} \}$,
the matrix of $A$ is
\[
A^{(\mathcal{B})} =
\left[ \begin{array}{cccccccc}
\gamma_1 &         &        &          &   &        &   & \beta_0 \\
1        & \gamma_2 &        &          &   &        &   & \beta_1 \\ 
         & 1       & \ddots &          &   &        &   & \vdots \\
         &         & \ddots & \gamma_J &   &        &   & \beta_{J-1} \\
         &         &        & 1        & 0 &        &   & \beta_J \\
         &         &        &          & 1 & \ddots &   & \vdots \\
         &         &        &          &   & \ddots & 0 & \beta_{n-2} \\
         &         &        &          &   &        & 1 & \beta_{n-1} \end{array} \right] .
\]
It is known that any $(n-1) \times (n-1)$ upper Hessenberg matrix with $1$'s on its
subdiagonal can be extended to an $n \times n$ upper Hessenberg matrix with $1$'s
on its subdiagonal that has any desired eigenvalues, by suitable choice of the 
last column; see, for instance, \cite{ParlettStrang}.  Thus we can choose $\beta_0 , \ldots ,
\beta_{n-1}$ so that $A^{(\mathcal{B})}$ has the given eigenvalues $\lambda_1 , \ldots , 
\lambda_n$.  If $B = [ b , u_1 , \ldots , u_{J-1} , v_1 , \ldots , v_{n-J} ]$, then
the matrix of $A$ in the standard basis is $B^{-1} A^{(\mathcal{B})} B$.
This completes the construction of the desired matrix $A$ and right-hand side vector $b$.

\section{Further Remarks}
We have introduced a method, Arnoldi-OR, that uses
$\max \{ \mbox{deg}(D(A)), \mbox{deg}(N(A)) \} -1$ extra steps of the Arnoldi algorithm
in order to find the optimal
(in $D(A )^{*} D(A)$-norm) approximation to $D(A )^{-1} N(A) b$ from the Krylov space
$\mbox{span} ( b, Ab, \ldots , A^{k-1} b )$.  To do this requires solving a $k + \max \{ \mbox{deg}
(D(A)), \mbox{deg}(N(A)) \}$ by $k$ least squares problem.  This can be solved via QR factorization with
Givens rotations, applying previous rotations to just the last column of the newly formed
matrix and choosing new rotations to eliminate entries below the diagonal in column $k$.
The residual norm in the least squares problem is equal to $\| N(A) b - D(A) Q_k y \|_2$,
so the approximate solution $Q_k y$ need not be formed until a desired tolerance is met.

For normal and near-normal matrices $A$, error bounds can be based on how well the
rational function $R(z) := D(z )^{-1} N(z)$ can be approximated by a polynomial of degree at 
most $k-1$ on the eigenvalues of $A$, independent of the eigenvalues of the upper
Hessenberg matrix that happens to be generated by the Arnoldi algorithm.
For highly nonnormal problems, {\em a priori} error bounds can be based on how well
$R(z)$ can be approximated by a polynomial of degree at most $k-1$ on
the numerical range $W(A)$ or on other $K$-spectral sets.

In many of the numerical experiments, one sees that there is not much difference 
between the Arnoldi-OR and the Arnoldi-FA convergence curves.  In fact, when used
for solving linear systems, a precise relationship between the two methods is 
known \cite{Brown,CullumG}:
\[
\| r_k^{FA} \|_2 = \frac{\| r_k^{OR} \|_2}{\sqrt{1 - ( \| r_k^{OR} \|_2 / 
\| r_{k-1}^{OR} \|_2 )^2}} ,
\]
where $r_k^{FA}$ denotes the residual in the Arnoldi-FA, or, FOM
approximation and $r_k^{OR}$ denotes the residual in the Arnoldi-OR, or,
GMRES approximation.  This precise relationship does not generally hold 
for rational functions of higher degree, but it would not be surprising if a similar
type of relation could be derived.  In this case, the optimality of Arnoldi-OR
could then be used as an aid in developing {\em a priori} error bounds for
Arnoldi-FA.


\begin{thebibliography}{99}
\bibitem{RKToolkit} {\sc M.~Berljafa, S.~Elsworth, and S.~G\"{u}ttel}, {\em Rational Krylov toolbox for MATLAB},
\verb+http://guettel.com/rktoolbox/+ .

\bibitem{Brown} {\sc P.~Brown},
{\em A theoretical comparison of the Arnoldi and GMRES algorithms},
{SIAM J. Sci. Stat. Comput.}, 20 (1991), pp.~58-78.

\bibitem{CGMM1} {\sc T.~Chen, A.~Greenbaum, C.~Musco, and C.~Musco},
{\em Error bounds for Lanczos-based matrix function approximation}, {ArXiv https://arxiv.org/abs/2106.09806},
{to appear in SIAM J. Matrix Anal. Appl.}.

\bibitem{CGMM2} {\sc T.~Chen, A.~Greenbaum, C.~Musco, and C.~Musco},
{\em Low-memory Krylov subspace methods for optimal rational matrix function 
approximation}, 
{SIAM J. Matrix Anal. Appl.}, 44 (2023), pp.~670-692.

\bibitem{CG} {\sc M.~Crouzeix, A.~Greenbaum},
{\em Spectral sets:  numerical range and beyond}, {SIAM J. Matrix Anal. Appl.},
40 (2019), pp.~1087-1101.

\bibitem{CP} {\sc M.~Crouzeix, C.~Palencia},
{\em The numerical range is a $(1{+}\sqrt2)$-spectral set}, {SIAM J. Matrix Anal. Appl.}, 
38 (2017), pp.~649-655.

\bibitem{CullumG} {\sc J.~Cullum, A.~Greenbaum},
{\em Relations between Galerkin and norm-minimizing iterative methods for solving linear 
systems}, {SIAM J. Matrix Anal. Appl.}, 17 (1996), pp.~223-247.

\bibitem{Driscoll} {\sc T.~Driscoll and L.~Trefethen},
{\em Schwarz-Christoffel Mapping}, {Cambridge}, 2002.

\bibitem{Embree} {\sc M.~Embree}, {\em How descriptive are GMRES convergence
bounds}, arXiv2209.01231v1, Sept. 2, 2022.

\bibitem{Gaier} {\sc D.~Gaier}, {\em Lectures on Complex Approximation}, {Birkh\"{a}user},
1980.

\bibitem{GPS} {\sc A.~Greenbaum, V.~Ptak, and Z.~Strako\v{s}},
{\em Any nonincreasing convergence curve is possible for GMRES}, {SIAM J.~Matrix Anal.~Appl.},
17 (1996), pp.~465-469.

\bibitem{GW} {\sc A.~Greenbaum and N.~Wellen},
{\em $K$-spectral sets}, arXiv2302.05535, submitted to SIAM J.~Matrix Anal.~Appl..

\bibitem{Moler} {\sc C.~Moler and C,~Van Loan},
{\em Nineteen dubious ways to compute the exponential of a matrix, twenty-five
years later}, {SIAM Review} 45 (2003), pp.~3-49.

\bibitem{Napoli} {\sc E.~Di Napoli, E.~Polizzi, and Y.~Saad},
{\em Efficient estimation of eigenvalue counts in an interval}, 
{Num.~Lin.~Alg.~Appl.}, 23 (2016), pp.~674-692.

\bibitem{ParlettStrang} {\sc B.~Parlett and N.~Strang},
{\em Matrices with prescribed Ritz values}, {Lin.~Alg.~Appl.}, 428 (2008), pp.~1725-1739.

\bibitem{eigtool} {\sc T. M. ~Wright, M. ~Embree}
{\em EigTool: a graphical tool for nonsymmetric eigenproblems}, {\texttt{http://www.cs.ox.ac.uk/pseudospectra/eigtool/}}

\end{thebibliography}
\end{document}